\magnification \magstep1
\outer\def\thm#1#2\par{\medbreak\noindent{\bf(#1)\enspace}\ignorespaces
{\sl#2}\ifdim\lastskip<\medskipamount\removelastskip\penalty55\medskip\fi}
\def\cond#1#2\par{\smallbreak\noindent\rlap{\rm(#1)}\ignorespaces
\hangindent=36pt\hskip36pt{\rm#2}\smallskip}
\def\claim#1#2\par{\smallbreak\noindent\rlap{\rm(#1)}\ignorespaces
\hangindent=30pt\hskip30pt{\sl#2}\smallskip}
\def\dfn#1{{\sl #1}}

\def\proof{\smallbreak\noindent{\sl Proof. }}



\def\qed{\hfill$\square$\medskip}
\def\sqr#1#2{{\vcenter{\vbox{\hrule height.#2pt
\hbox{\vrule width.#2pt height #1pt \kern#1pt
\vrule width.#2pt}
\hrule height.#2pt}}}}
\def\square{\mathchoice\sqr56\sqr56\sqr{2.1}3\sqr{1.5}3}

\def\ref#1#2{\item{#1.}#2}
\outer\def\beginsection#1\par{\vskip0pt plus.3\vsize
   \vskip0pt plus-.3\vsize\bigskip\bigskip\vskip\parskip
   \message{#1}\centerline{\bf#1}\nobreak\smallskip\noindent}

\def\st{\,|\,}
\def\junk#1{}

\font\smallrm=cmr8

\def\Ccal{{\cal C}}
\def\Mcal{{\cal M}}

\def\reffc{1}

\nopagenumbers
\phantom{a}
\bigskip
\centerline{\bf REDUCIBILITY IN THE FOUR-COLOR THEOREM}
\bigskip\bigskip
\baselineskip=11pt

\centerline{Neil Robertson$^{*1}$\vfootnote{$^*$}{\smallrm
Research partially 
performed under a consulting agreement with Bellcore, and partially
supported by DIMACS Center, 
Rutgers University, New Brunswick, New Jersey  08903, USA.}
\vfootnote{$^1$}{\smallrm Partially supported
by NSF under Grant No. DMS-8903132 and by ONR under Grant No. 
N00014-92-J-1965.
}}
\centerline{Department of Mathematics}
\centerline{Ohio State University}
\centerline{231 W. 18th Ave.}
\centerline{Columbus, Ohio  43210, USA}
\bigskip

\centerline{Daniel P. Sanders$^{2}$\vfootnote{$^2$}{\smallrm
Partially supported by DIMACS and by ONR under Grant No. 
N00014-93-1-0325.
}}
\centerline{School of Mathematics}
\centerline{Georgia Institute of Technology}
\centerline{Atlanta, Georgia  30332, USA}
\bigskip

\centerline{Paul Seymour}
\centerline{Bellcore}
\centerline{445 South St.}
\centerline{Morristown, New Jersey  07960, USA}
\bigskip

\centerline{and}
\bigskip

\centerline{Robin Thomas$^{*3}$\vfootnote{$^3$}{\smallrm
Partially supported
by NSF under Grant No. DMS-9303761 and by ONR under Grant No. 
N00014-93-1-0325.
}}
\centerline{School of Mathematics}
\centerline{Georgia Institute of Technology}
\centerline{Atlanta, Georgia  30332, USA}
\bigskip\bigskip\bigskip

\centerline{\bf ABSTRACT}
\bigskip

In [{\it J.~Combin.\ Theory Ser.\ B} 70 (1997), 2-44] we gave a simplified proof
of the Four-Color Theorem. The proof is computer-assisted in the sense
that for two lemmas in the article we did not give proofs, and instead
asserted that we have verified those statements using a computer. Here
we give additional details for one of those lemmas, and we include the
original computer programs and data as ``ancillary files" accompanying this
submission. 


\vfill
\noindent 29 March 1995.
 Revised 27 January 1997
\vfil
\eject

\footline{\hss\tenrm\folio\hss}
\baselineskip 22pt

\beginsection 1. CONFIGURATIONS

We assume familiarity with [\reffc].
The purpose of this manuscript
is to provide more details about the proof of [\reffc, theorem (3.2)].
As a first step we need to explain how configurations are stored.

Let $R$ be a circuit with vertices $1,2,\ldots,r$ in order, let
$e_1,e_2,\ldots,e_r$ be the edges of $R$, in order and such that
for $i=2,3,\ldots,r$, the ends of $e_i$ are $i$ and $i-1$. Let
$\kappa,\lambda:E(R)\to \{-1,0,1\}$ be two colorings of $R$. We
say that $\kappa$ and $\lambda$ are
\dfn{similar} if $\{\kappa^{-1}(-1),\kappa^{-1}(0),\kappa^{-1}(1)\}=
\{\lambda^{-1}(-1),\lambda^{-1}(0),\lambda^{-1}(1)\}$. We say
that $\kappa$ is \dfn{canonical} if either $\kappa(e)=0$ for every
edge $e$ of $R$, or there exists an integer $k$ such that $1\le k<r$,
$\kappa(e_r)=\kappa(e_{r-1})=\cdots=\kappa(e_{k+1})=0$ and
$\kappa(e_k)=1$. Clearly every coloring of $R$ is similar to
a unique canonical coloring.
The \dfn{code} of a coloring $\kappa:E(R)\to \{-1,0,1\}$ is
$\sum_{i=1}^r\kappa'(e_i)3^{i-1}$, where $\kappa'$ is the canonical
coloring similar to $\kappa$. 

Let $K$ be a configuration, let $G$ be the free completion of $K$
with ring $R$, and let either $X$ be empty or
a contract for $K$. Let the vertices of $G$ be $1,2,\ldots,n$, where
$1,2,\ldots,r$ are the vertices of $R$ in clockwise
order around the infinite region of $G$. Let $\Ccal(K)$ be the set
of all restrictions to $E(R)$ of tri-colorings of $G$, and let
$\Ccal'(K)$ be the maximal consistent subset of $\Ccal^*-\Ccal(K)$,
where $\Ccal^*$ is the set of all mappings of $E(R)\to\{-1,0,1\}$.
A matrix $A=(a_{i,j})$ is
a \dfn{configuration matrix of $K$} if
\item{(i)}for $i=1,2,\ldots,n$, $a_{i,0}$ is the degree of vertex
$i$ of $G$, and $a_{i,1},a_{i,2},\ldots,a_{i,a_{i,0}}$ are the neighbors
of $i$ listed in clockwise order as they appear around $i$; moreover,
if $i\le r$ then $a_{i,1}$ and $a_{i,a_{i,0}}$ belong to $V(R)$,
\item{(ii)}$a_{0,0}=n$ and $a_{0,1}=r$,
\item{(iii)}$a_{0,2}$ and $a_{0,3}$ are the numbers of canonical 
colorings in $\Ccal(K)$ and $\Ccal'(K)$, respectively,
\item{(iv)}$a_{0,4}=|X|$, and $X=\{f_1,f_2,\ldots,f_k\}$, where
$k=|X|$ and for $i=1,2,\ldots,k$, $f_i$ has ends $a_{0,2i+3}$
and $a_{0,2i+4}$.

\noindent Each of the $633$ good configurations is presented in terms
of a configuration matrix, and so we need to verify that an 
input matrix is indeed a configuration matrix of some configuration.
Let $r$ and $n$ be integers, let $A=(a_{i,j})$ be an integer matrix
with rows corresponding to $i=0,1,\ldots,n$,
and let us consider the following conditions, where for notational
convenience we put $d_i=a_{i,0}$.
\item{(1)}$2\le r<n$,
\item{(2)}$3\le d_i\le n-1$ for all $i=1,2,\ldots,r$, and 
$5\le d_i\le n-1$ for all $i=r+1,r+2,\ldots,n$.
\item{(3)}$1\le a_{i,j}\le n$ for $i=1,2,\ldots,n$ and $j=1,2,\ldots,
d_i$,
\item{(4)}if $i=1,2,\ldots,r$, then $a_{i,1}=i+1$ (or $1$ if $i=r$),
$a_{i,d_i}=i-1$ (or $r$ if $i=1$), and
$r+1\le a_{i,j}\le n$ for 
$j=2,3\ldots,d_i-1$,
\item{(5)}$d_1+d_2+\cdots+d_n=6n-6-2r$,
\item{(6)}for every $i=r+1,r+2,\ldots,n$ there exist at most two
integers $j$ such that $a_{i,j}> r$ and $a_{i,j+1}\le r$, and if
there are two then $a_{i,j+2}>r$ for both such integers (where
$a_{i,d_i+1}$ and $a_{i,d_i+2}$ mean $a_{i,1}$ and $a_{i,2}$, 
respectively), and
\item{(7)}let $i=1,2,\ldots,n$, let $j=1,2,\ldots,d_i-1$
if $i\le r$ and $j=1,2,\ldots,d_i$ if $i>r$, and let $k=a_{i,j}$.
Then there exists an integer $p$ such that $a_{i,j+1}=a_{k,p}$
(or $a_{i,1}=a_{k,p}$ if $j=d_i$), and
$i=a_{k,p+1}$ (or $i=a_{k,1}$ if $p=d_k$, in which case $k>r$).

\thm{1.1}Let $r,n$ and $A=(a_{i,j})$ satisfy (1)--(7). Then there
exist a configuration $K$ and a free completion $G$ of $K$ with
ring $R$ such that (i) and (ii) hold.

\proof It is straightforward to construct a graph $G$ with
vertex-set $\{1,2,\ldots,n\}$ such that for $i=1,2,\ldots,n$,
the neighbors of $i$ are $a_{i,1},a_{i,2},\ldots,a_{i,d_i}$.
From (7) we deduce that the cyclic orderings 
$a_{i,1},a_{i,2},\ldots,a_{i,d_i}$ of the neighbors of $i$
for $i=1,2,\ldots,n$ define an embedding of $G$ into a surface
$\Sigma$ such that every face is a triangle, except one bounded
by a circuit $R$. By (4) the vertices of $R$ are $1,2,\ldots,r$
in order. From (5) we deduce by Euler's formula that $\Sigma$
is the sphere, and so $G$ may be regarded as a near-triangulation.
Let $K$ be defined by 
$G(K)=G\backslash V(C)$, and for $v\in V(K)$ let $\gamma_K(v)=
d_G(v)$. We claim that $K$ is a configuration. By (4)
$R$ is an induced circuit of $G$, and hence $G(K)$ is connected.
Thus $G(K)$ is a near-triangulation, and we must verify conditions
(i), (ii) and (iii) in the definition of a configuration. Condition
(i) follows from (6), condition (ii) follows from (2), because
every vertex on the infinite region of $G(K)$ is adjacent to a vertex
of $R$, and (iii) follows from (1), because $r$ is the ring-size
of $K$. \qed

\beginsection 2. EXTENDABLE COLORINGS

The objective of this section is to explain how we compute $\Ccal(K)$.
We compute all tri-colorings of $G$ and record their restrictions
to $E(R)$. We first number the edges of $G$ as $e_1,e_2,\ldots,e_m$,
where $m=3(n-1)-r$ and the edges of $R$ are $e_1,e_2,\ldots,e_r$.
All that matters for the correctness is that $e_m$ and $e_{m-1}$
are on a common triangle. Using the algorithm below we compute 
all mappings $c:E(G)-E(R)\to
\{1,2,4\}$ such that $c(e_m)=1$, $c(e_{m-1})=2$ and
$c(e)\not=c(f)$ if $e$ and $f$ are on a
common triangle. During the course of
the algorithm we maintain a variable $F_i$ defined for $i<m-1$ as the
set of all $c(e_j)$ such that $j>i$ and $e_i$ and $e_j$ are on a common
triangle. At the begining we set $c(e_m)=1$, $c(e_{m-1})=2$,
$F_{m-1}=\{1,4\}$ and $j=m-1$, and keep repeating steps 1, 2,
3 below.

\noindent {\bf Step 1.} While $c(e_j)\in F_j$ we keep repeating
steps (i) and (ii) below.
\item{(i)}Double $c(e_j)$, and
\item{(ii)}while $c(e_j)=8$ repeat the following steps:
\itemitem{(a)}if $j\ge m-1$ terminate computation,
\itemitem{(b)}increase $j$ by one and double $c(e_j)$.

\noindent {\bf Step 2.} If $j=r+1$ then a tri-coloring of $G$ can
be read off from $c$. Record the code of its restriction to $E(R)$.
Double $c(e_j)$ and while $c(e_j)=8$ repeat steps (a) and (b)
above.

\noindent {\bf Step 3.} If $j>r+1$ decrease $j$ by one, set $c(e_j)=1$,
and compute $F_j$.

We record the codes of restrictions of tri-colorings of $G$ to $E(R)$
in an array called ``live", so that for $i=0,1,\ldots,(3^{r-1}-1)/2$,
live$[i]=0$ if some (and hence every) coloring of $E(R)$ with code
$i$ is the restriction to $E(R)$ of a tri-coloring of $G$, and
live$[i]=1$ otherwise.

\beginsection 3. CONSISTENT SETS

We now explain how we compute $\Ccal'(K)$. Let $K,G,R,\Ccal^*$ be as in
Section 1. We say that a coloring $\kappa$ of $R$ is \dfn{balanced}
if $|\kappa^{-1}(-1)|$, $|\kappa^{-1}(0)|$, $|\kappa^{-1}(1)|$ 
and $r$ all have the same parity. It is easy to see that
every member of $\Ccal'(K)$ is balanced. Let
$M=\{(m_1,\mu_1),(m_2,\mu_2),\ldots,\allowbreak(m_k,\mu_k)\}$ 
be a signed
matching in $R$. We say that $M$ is \dfn{balanced} if
$r-\sum_{i=1}^k(\mu_i-1)/2$ is even.
Let $\Ccal_0=\Ccal^*-\Ccal(K)$, and let $\Mcal_0$ be
the set of all balanced signed matchings in $R$. Let $i\ge0$ be an
integer, and assume that $\Mcal_0,
\Mcal_1,\ldots,\Mcal_i$ and $\Ccal_0,\Ccal_1,\ldots,\Ccal_i$
have already been defined. Define $\Mcal_{i+1}$ to be the set
of all signed matchings $M\in\Mcal_i$ such that $\Ccal_i$ contains
every coloring
$\kappa$ of $R$ that $\theta$--fits $M$ for some $\theta\in\{-1,0,1\}$,
and let $\Ccal_{i+1}$ be the set
of all colorings $\kappa\in\Ccal_i$ such that for every $\theta
\in\{-1,0,1\}$ there is a signed matching $M\in\Mcal_{i+1}$
such that $\kappa$ $\theta$--fits $M$. We need the following
proposition.

\thm{3.1}If $\Ccal_i=\Ccal_{i+1}$, then $\Ccal'(K)=\Ccal_i$.

\proof We first show that $\Ccal_i$ is consistent. To this end
we notice that $\Mcal_{i+1}=\Mcal_{i+2}$. Let $\kappa\in\Ccal_i$, and
let $\theta\in\{-1,0,1\}$. Since $\kappa\in\Ccal_{i+1}$ we deduce
that there exists a signed matching $M\in\Mcal_{i+1}$ such that
$\kappa$ $\theta$--fits $M$. Since $M\in\Mcal_{i+2}$ by the
above observation, if a coloring $\kappa'$ $\theta$--fits $M$,
then $\kappa'\in\Ccal_{i+1}=\Ccal_i$, as desired.

To complete the proof we must show that $\Ccal'(K)\subseteq\Ccal_j$
for all $j=0,1,\ldots$. We proceed by induction. Clearly
$\Ccal'(K)\subseteq\Ccal_0$. Assume now that $\Ccal'(K)\subseteq
\Ccal_j$ for some integer $j\ge 0$; we wish to show that 
if $\kappa\in\Ccal_j-\Ccal_{j+1}$, then $\kappa\not\in\Ccal'(K)$.
Let $\kappa$ be as stated. Then there exists $\theta\in\{-1,0,1\}$
such that $\kappa$ $\theta$--fits no $M\in\Mcal_{j+1}$.
If $\kappa$ $\theta$--fits no signed matching in $R$ then $\kappa\not\in
\Ccal'(K)$, and so we may assume that $\kappa$ $\theta$--fits
a signed matching $M$. If $M\not\in\Mcal_{0}$, then $\kappa$
is unbalanced, and hence $\kappa\not\in\Ccal'(K)$. We may therefore
assume that $M\in\Mcal_{0}$. Then $M\in \Mcal_k-\Mcal_{k+1}$
for some integer $k$ with $0\le k\le j$. Hence there exists
$\theta'\in\{-1,0,1\}$ and a coloring $\kappa'\not\in\Ccal_k$
$\theta'$-fitting $M$. By replacing $\kappa'$ by a similar
coloring if necessary we may assume that $\theta=\theta'$.
Since $\kappa'\not\in\Ccal'(K)$ by the induction hypothesis,
we deduce that $\kappa\not\in\Ccal'(K)$, as desired. \qed

To compute $\Ccal'(K)$ we iteratively compute $\Mcal_i$ and
$\Ccal_i$ until $\Ccal_i=\Ccal_{i+1}$. Instead of $\Ccal_i$
we compute the codes of members of $\Ccal_i$, and store this
information by updating the array ``live". To complete the
description we need to explain how we store and compute $\Mcal_i$.
To this end we need the following definitions. Let
$M=\{(m_1,\mu_1),(m_2,\mu_2),\ldots,(m_k,\mu_k)\}$ be a signed
matching in $R$, where for $i=1,2,\ldots,k$, $m_i=\{a_i,b_i\}$,
$b_i<a_i$ and $a_1=\max\{a_1,a_2,\ldots,a_k\}$.
We define the \dfn{code} of $M$ to be
$\sum_{i=1}^k\left(3^{a_i-1}+\mu_i3^{b_i-1}\right)$ if $a_1<r$, and
$(3^r-1)/2-\sum_{i=1}^k\left(3^{a_i-1}+(3-\mu_i)3^{b_i-1}/2\right)$ 
otherwise. We define $(h_2,h_3,\ldots,h_k)$, the 
\dfn{choice sequence of $M$}, by $h_i=2(3^{a_i-1}+\mu_i3^{b_i-1})$ if
$a_1<r$ and $h_i=3^{a_i-1}+\mu_i3^{b_i-1}$ otherwise. The following
is straightforward to verify.

\thm{3.2}Let $M$ be a signed matching with code $c$ and choice
sequence $(h_2,h_3,\ldots,h_k)$. Then 
$\left\{\left|c+\sum_{i=2}^k\epsilon_i h_i\right|\,:\,
\epsilon_i\in\{0,1\}\right\}$ is
the set of codes of all colorings of $R$ that $\theta$--fit
$M$ for some $\theta\in\{-1,0,1\}$. Moreover, let $a_1$ be as above,
and let $\kappa$ be a canonical coloring $\theta$--fitting $M$
for some $\theta\in\{-1,0,1\}$. Let the code of $\kappa$ be
$d=c+\sum_{i=2}^k\epsilon_i h_i$, where $\epsilon_i\in\{0,1\}$
for $i=2,3,\ldots,k$. If $a_1<r$ then $\theta=0$. If $a_1=r$ then
$\theta=1$ if $d<0$, and $\theta=-1$ if $d>0$.

The function ``augment" generates all the members of $\Mcal_0$
always in the same order, say $M_1,M_2,\ldots,M_p$. We use the 
bits of the array ``real" to store $\Mcal_i$; that is, after the
$i$th iteration the $j$th bit of ``real" is $1$ if and only if
$M_j\in\Mcal_i$. To update ``real" we run through all bits of
``real" that are currently set to $1$, generate all the codes of
colorings as in (3.2), and if for some of them the corresponding
entry in ``live" is zero, we set the current bit of ``real" to
zero. Also, if none of the corresponding entries of ``live" are
zero, we mark each such entry by $\theta$, where $\theta$ is
as in the second half of (3.2). To update ``live"
(that is, to compute $\Ccal_{i+1}$ from $\Ccal_i$) we run through
all nonzero entries of ``live" and set to zero all those
that were not marked by every $\theta\in\{-1,0,1\}$ (except
live$[0]$, which is exceptional).

\beginsection 4. CONTRACTS

It remains to explain how we verify that a proposed contract
is indeed a contract. Let $K,G,R$ be as in Section 1,
and let $X$ be a contract for $K$ as specified by condition
(iv) in the definition of a configuration matrix. Most of
the conditions in the definition of a contract are straightforward
to verify, and so we only explain how we
check that no coloring in $\Ccal'(K)$ is the restriction
to $E(R)$ of a tri-coloring of $G$ modulo $X$.  We do that by 
computing
all tri-colorings of $G$ modulo $X$ using an algorithm similar
to the one described in Section 2. More precisely, the algorithm
runs as follows.
Let $e_1,e_2,\ldots,e_m$ be the numbering of $E(G)$ described
in Section 2. We compute a mapping $c:E(G)-X\to
\{1,2,4\}$ such that $c'$ defined by $c'(e)=\lfloor c(e)/2\rfloor-1$
is a tri-coloring of $G$ modulo $X$. For an integer $i$ with 
$1\le i\le m$ and $e_i\not\in X$
let $D_i$ be the set of all $e_j$, where $i<j\le m$ and $e_i$,
$e_j$ belong to a triangle none of whose edges belong to $X$.
Let $S_i$ be the set of all $e_j$, where $i<j\le m$, $e_j\not\in X$,
and $e_i$, $e_j$ belong to a triangle whose third edge belongs to $X$.
Let $s$ be the maximum integer with $s\le m$ and $e_{s}\not\in X$,
and let $s'$ be the maximum integer with $s'<s$ and $e_{s'}\not\in X$.
During the course of
the algorithm we maintain a variable $F_i$ defined for $i<s'$ as 
$\{c(f)\st f\in D_i\}\cup\bigcup\left\{\{1,2,4\}-\{c(f)\}\st f\in S_i
\right\}$.
At the beginning we set $c(e_s)=1$, $c(e_{s'})=1$,
$F_{s'}=\{4\}\cup\{c(f)\st f\in D_{s'}\}\cup\bigcup\left\{\{1,2,4\}-
\{c(f)\}\st
f\in S_{s'}\right\}$ and $j=s-1$, and keep repeating steps 1, 2,
3 below.

\noindent {\bf Step 1.} While $c(e_j)\in F_j$ we keep repeating
steps (i) and (ii) below.
\item{(i)}Double $c(e_j)$, and
\item{(ii)}while $c(e_j)=8$ repeat the following steps:
\itemitem{(a)}set $j$ to the smallest $j'$ with $s\ge j'>j$ and 
$e_{j'}\not\in X$ (or $s+1$ if no such $j'$ exists),
\itemitem{(b)}if $j\ge s$ terminate computation; otherwise double $c(e_j)$.

\noindent {\bf Step 2.} If $j=1$ then $c$ can be converted to
a tri-coloring of $G$ modulo $X$. Verify that its restriction to $E(R)$
does not belong to $\Ccal'(K)$. Double $c(e_j)$ and while $c(e_j)=8$
repeat steps (a) and (b) above.

\noindent {\bf Step 3.} If $j>1$ decrease $j$ by one, set $c(e_j)=1$, 
and compute $F_j$.

\beginsection ACKNOWLEDGMENT

We would like to express our thanks to Tom Fowler for reading this manuscript
and the corresponding program.

\beginsection REFERENCE


\item{\reffc.}N. Robertson, D. P. Sanders, P. D. Seymour and
R. Thomas, The Four-Colour Theorem, 
to appear in {\it J.\ Combin.\ Theory Ser.\ B}.

\end